\documentclass{svjour3}                     
\smartqed  
\usepackage{graphicx}
 \usepackage{mathptmx}      
\usepackage{amsmath,amsfonts,amssymb}
\usepackage{latexsym}
\begin{document}

\title{Generalized Fourier transform on Ch\'{e}bli-Trim\`{e}che hypergroups\thanks{ } }


\author{Chokri Abdelkefi \and Abdessattar Jemai
}


\institute{Chokri Abdelkefi
              \\Department of Mathematics, Preparatory
Institute of Engineer Studies of Tunis,\\ 1089 Monfleury Tunis,
Tunisia\\\email{chokri.abdelkefi@ipeit.rnu.tn}\\ \\
           Abdessattar Jemai \at
              Department of Mathematics, Faculty of Sciences
of Tunis, \\1060 Tunis, Tunisia \\
              \email{jemai\_abdessattar@yahoo.fr}  }

\date{Received: date / Accepted: date}

\maketitle

\begin{abstract}
In this paper, we prove the Hardy-Littlewood inequality for the
generalized Fourier transform on Ch\'{e}bli-Trim\`{e}che hypergroups
and we study in the particular case of the Jacobi hypergroup the
integrability of this transform on Besov-type spaces.
\end{abstract}

\section{Introduction}
\label{intro}  We consider the Ch\'{e}bli-Trim\`{e}che hypergroup
$(\mathbb{R}_{+}, \ast(A))$ associated with the function $A$ which
depends on a real parameter $\alpha>-\frac{1}{2}$ (see next
section). We prove the Hardy-Littlewood inequality for the
generalized Fourier transform $\mathcal{F}(f)$ of a function $f$ in
$L^{p}(\mathbb{R}_{+}, A(x)dx)$, $1<p\leq2$. Next, inspired by the
definition of usual Besov spaces and Besov-Dunkl spaces (see [2,
5]), we define the Besov-type spaces for Ch\'{e}bli-Trim\`{e}che
hypergroup denoted by $\mathcal{B}_{\gamma,\alpha}^{p,q}$, as the
subspace of functions $f\in L^{p} (\mathbb{R}_{+},A(x)dx)$
satisfying
\begin{eqnarray*} \int_{0}^{+\infty}\Big(\frac{\omega_{A,p}(f,x)}{x^{\gamma}}\Big)^{q}\frac{dx}{x}<+\infty \quad\mbox{if } q<+\infty \end{eqnarray*}
and
\begin{eqnarray*}\sup_{x\in]0,+\infty[}\frac{\omega_{A,p}(f,x)}{x^{\gamma}}<+\infty \quad \mbox{if } q=+\infty,\end{eqnarray*}
where $\omega_{A,p}(f,x)=\displaystyle{ \|\tau_{x}(f)-f\|_{A,p}}$ is
the modulus of continuity of first order of $f$ with $\tau_{x}$ the
generalized translation operators, $x\in \mathbb{R}_+$ (see next
section). We establish in the particular case of Jacobi hypergroup
further results concerning integrability of the generalized Fourier
transform $\mathcal{F}(f)$ of a function $f$ when $f$ belongs to a
suitable Besov-type spaces. Analogous results have been obtained for
the theory of Dunkl operators in [1, 3, 4].

The contents of this paper are as follows.\\
 In section 2, we collect some results about harmonic analysis
 on Ch\'{e}bli-Trim\`{e}che hypergroups.\\ In section 3, we prove the Hardy-Littlewood inequality for the
generalized Fourier transform on Ch\'{e}bli-Trim\`{e}che hypergroups
and we study in the particular case of the Jacobi hypergroup the
integrability of this transform on Besov-type spaces.
 \par Along this paper we use $c$ to
denote a suitable positive constant which is not necessarily the
same in each occurrence. Furthermore, we denote by

$\bullet\quad \mathbb{C}_{\ast,c}(\mathbb{R})$ the space of even
continuous functions on $\mathbb{R}$, with compact support.

$\bullet\quad \mathcal{D}_\ast(\mathbb{R})$ the space of even
$C^{\infty}$-functions on $\mathbb{R}$
 with compact support.

\section{Preliminaries}
\label{sec:1}  In this section, we recall some notations and results
about
harmonic analysis on Ch\'{e}bli-Trim\`{e}che hypergroups and we refer for more details to the articles [6, 9, 11, 12].\\

Let $A$ be the Ch\'{e}bli-Trim\`{e}che function defined on
$\mathbb{R}_{+}$ and satisfying the following conditions.
\begin{itemize}
\item [i)] $A(x)=x^{2\alpha+1}B(x)$, with $\alpha>-\frac{1}{2}$, and $B$ an even $C^{\infty}$-function on $\mathbb{R}$ such that $B(x)\geq1$
for all $x\in \mathbb{R}_{+}$.
\item [ii)] $A$ is increasing and unbounded.
\item [iii)] $\displaystyle\frac{A'}{A}$ is decreasing on $\mathbb{R}_{+}^{*}=]0,+\infty[$ and
$\displaystyle\lim_{x\mapsto+\infty}\frac{A'(x)}{A(x)}=2\rho\geq0$.
\item [iv)] There exists a constant $\eta>0$ such that for all $x\in [x_{0}, +\infty[,\, x_{0}>0$, we have
\[
\frac{A'(x)}{A(x)}= \left\{
\begin{array}{ll}
2\rho+ e^{-\eta x }F(x)& \mbox{, if } \rho>0 \\
\frac{2\alpha+1}{x}+e^{-\eta x }F(x)& \mbox{, if } \rho=0,
\end{array}
\right.\] where $F$ is a $C^{\infty}$-function bounded together with
its derivatives.
\end{itemize}

We consider the Ch\'{e}bli-Trim\`{e}che hypergroup $(\mathbb{R}_{+},
\ast(A))$ associated with the function $A$. We note that it is
commutative with neutral element $0$ and the identity mapping is the
involution. The Haar measure $m$ on $(\mathbb{R}_{+}, \ast(A))$ is
absolutely continuous with respect to the Lebesgue measure and can
be choosen to have the Lebesgue density $A$.
\begin{remark}
If $A(x)=2^{2\rho}(\sinh x)^{2\alpha+1}(\cosh x)^{2\beta+1}$, with
$\alpha\geq \beta\geq-\frac{1}{2}$, $\alpha\neq-\frac{1}{2}$ and
$\rho=\alpha+\beta+1$, $(\mathbb{R}_{+}, \ast(A))$ is called the
Jacobi hypergroup.
\end{remark}

Let $\Delta$ be the differential operator on $\mathbb{R}^{*}_{+}$
given by
$$\Delta=\frac{d^{2}}{dx^{2}}+\frac{A'(x)}{A(x)}\frac{d}{dx}.$$
The solution $\varphi_{\lambda}, \lambda\in\mathbb{C}$, of the
differential equation
$$\left\{
\begin{array}{ll}
\Delta u(x)= -(\lambda^{2}+\rho^{2})u(x), \\
u(0)=1,\frac{d}{dx}u(0)=0,
\end{array}
\right.$$ is multiplicative on $(\mathbb{R}_{+}, \ast(A))$ in the
sense that
$$\forall x,y\in \mathbb{R}_{+}, \, \int_{\mathbb{R}_{+}} \varphi_{\lambda}(t)
\;d(\delta_{x}\ast\delta_{y})(t)=\varphi_{\lambda}(x)\varphi_{\lambda}(y),$$
where $\delta_{x}$ is the point mass at $x$ and
$\delta_{x}\ast\delta_{y}$ is a probability measure which is
absolutely continuous with respect to the measure $m$ and satisfies
$$\mbox{supp}\, \delta_{x}\ast\delta_{y}=[|x-y|, x+y].$$

We list some known properties of the characters $\varphi_\lambda$ of
the hypergroups.
\begin{itemize}
\item [i)] For each $\lambda\in \mathbb{C}$, the function $x\mapsto\varphi_{\lambda}(x)$ is an even
$C^{\infty}$-function on $\mathbb{R}$ and  for each $x\in
\mathbb{R}_{+}$, the fonction $\lambda\mapsto\varphi_{\lambda}(x)$
is an  entire function on $\mathbb{C}$.
\item [ii)] For every $\lambda\in \mathbb{C}$, the function $\varphi_{\lambda}$ admits the integral representation
$$\forall x\in\mathbb{R}_{+}^{\ast}, \quad \varphi_{\lambda}(x)=\int_{0}^{x} K(x,y) \cos(\lambda y)dy.$$
Where $K(x,.)$ is a positive even $\mathbb{C}^{\infty}$-function  on
$]-x, x[$ with  support in $[-x, x]$.
\end{itemize}

\begin{remark}
In the Jacobi hypergroup (see Remark 1), we have for all
$x\in\mathbb{R}_{+}$ and $\lambda\in\mathbb{C}$,
$$\varphi_{\lambda}(x)=\varphi_{\lambda}^{(\alpha,\beta)}(x)=
{}_{\phantom{1}2}\textrm{F}_{1}(\frac{1}{2}(\rho-i\lambda),\frac{1}{2}(\rho+i\lambda),\alpha+1;-\sinh^{2}x),$$
where ${}_{\phantom{1}2}\textrm{F}_{1}$ is the Gauss hypergeometric
function (see [9]). The function
$\varphi_{\lambda}^{(\alpha,\beta)}(x)$ is the Jacobi function and
it satisfies for all $ \lambda\in\mathbb{R}$ and
$t\in\mathbb{R}_{+}^{*}$
\begin{eqnarray}
|1- \varphi_{\lambda}^{(\alpha,\beta)}(t)|\geq c\, \min \{1,
(\lambda t)^{2}\},
\end{eqnarray}
where $c$ is constant which depends only on $\alpha$ and $\beta$
(see [7, 8]).
\end{remark}

For every  $p\in[1, +\infty]$, we denote by
$L^{p}_{A}(\mathbb{R}_{+})$ the space $L^{p}(\mathbb{R}_{+},
A(x)dx)$ and by $L^{p}_{c}(\mathbb{R}_{+})$ the space
$L^{p}(\mathbb{R}_{+}, \frac{d\lambda}{|c(\lambda)|^{2}})$ where
$|c(\lambda)|^{-2}$ is an even  continuous function on $\mathbb{R}$,
satisfying the estimates: There exist positive constants $k, k_{1},
k_{2}$ such that
\begin{itemize}
\item [i)] If $\rho=0$ and $\alpha>0$ then
\begin{eqnarray}
k_{1}|\lambda|^{2\alpha+1}\leq |c(\lambda)|^{-2}\leq
k_{2}|\lambda|^{2\alpha+1}, \quad \lambda\in\mathbb{C}.
\end{eqnarray}
\item [ii)] If $\rho>0$ and $\alpha>-\frac{1}{2}$ then
\begin{eqnarray}
k_{1}|\lambda|^{2\alpha+1}\leq |c(\lambda)|^{-2}\leq
k_{2}|\lambda|^{2\alpha+1}, \quad \lambda\in\mathbb{C},\;
|\lambda|>k,
\end{eqnarray}
and
\begin{eqnarray}
k_{1}|\lambda|^{2}\leq |c(\lambda)|^{-2}\leq k_{2}|\lambda|^{2},
\quad \lambda\in\mathbb{C}, \;|\lambda|\leq k.
\end{eqnarray}
\end{itemize}
We use $\|.\|_{A,p}$ and $\|.\|_{c,p}$ as a shorthand respectively
of $\|.\|_{L^{p}_{A}(\mathbb{R}_{+})}$ and
$\|.\|_{L^{p}_{c}(\mathbb{R}_{+})}$.\\

For $f\in L^{1}_{A}(\mathbb{R}_{+})$ the generalized Fourier
transform of $f$ is given by
$$\mathcal{F}(f)(\lambda)=\int_{\mathbb{R}_{+}}f(x)\varphi_{\lambda}(x)A(x)dx.$$
The generalized Fourier transform satisfies the following
properties.
\begin{itemize}
\item [i)] For $f\in L^{1}_{A}(\mathbb{R}_{+})$, we have
\begin{eqnarray}
\|\mathcal{F}(f)\|_{c,\infty}\leq\|f\|_{A,1}
\end{eqnarray}
\item [ii)] For $f$ in $ L^{1}_{A}(\mathbb{R}_{+})$ such that $\mathcal{F}(f)$ belongs to $ L^{1}_{c}(\mathbb{R}_{+})$, we have
the following inversion formula for the transform $\mathcal{F}$
\begin{eqnarray*}
f(x)=\int_{\mathbb{R}_{+}}\mathcal{F}(f)(\lambda)\varphi_{\lambda}(x)\frac{d\lambda}{|c(\lambda)|^{2}},\;
a.e.
\end{eqnarray*}
\item [iii)] (Plancherel formula) For all $f\in\mathcal{D}_{\ast}(\mathbb{R})$, we have
\begin{eqnarray}
\int_{\mathbb{R}_{+}}|f(x)|^{2}A(x)dx=\int_{\mathbb{R}_{+}}
|\mathcal{F}(\lambda)|^{2}\frac{d\lambda}{|c(\lambda)|^{2}}.
\end{eqnarray}
The transform $\mathcal{F}$ can be uniquely extended to an isometric
isomorphism from $ L^{2}_{A}(\mathbb{R}_{+})$ onto $
L^{2}_{c}(\mathbb{R}_{+})$.
\end{itemize}
For $1\leq p\leq2$, we denote by $p'$ the conjugate of $p$. From
(2.5), (2.6) and the Marcinkiewicz interpolation theorem (see [10]),
we obtain for $f\in L^{p}_{A}(\mathbb{R}_{+})$
\begin{eqnarray}
\mathcal{F}(f)\in L^{p'}_{c}(\mathbb{R}_{+}).
\end{eqnarray}

For $ x\in \mathbb{R}_{+}$ and
$f\in\mathbb{C}_{\ast,c}(\mathbb{R})$, the generalized $x$-translate
of $f$ is defined by
$$\forall y\in \mathbb{R}_{+}, \quad \tau_{x}f(y)=\int_{\mathbb{R}_{+}}f(t)d(\delta_{x}\ast\delta_{y})(t),$$
and we have $  \tau_{x}f(0)=f(x).$

The generalized translation operators $\tau_{x}$, $ x\in
\mathbb{R}_{+}$, satisfy the following properties.
\begin{itemize}
\item [i)] For all $x,y\in\mathbb{R}_{+}$ and $\lambda\in\mathbb{C}$, we have the product formula
\begin{eqnarray*}
\tau_{x}\varphi_{\lambda}(y)=\varphi_{\lambda}(x)\varphi_{\lambda}(y).
\end{eqnarray*}
\item [ii)] For $f\in\mathcal{D}_{\ast}(\mathbb{R})$ and $x\in\mathbb{R}_{+}$, the function $y\mapsto\tau_{x}f(y)$ belongs to
$\mathcal{D}_{\ast}(\mathbb{R})$ and we have
\begin{eqnarray}
\forall \lambda\in\mathbb{R}_{+}, \quad
\mathcal{F}(\tau_{x}f)(\lambda)=\varphi_{\lambda}(x)\mathcal{F}f(\lambda).
\end{eqnarray}
\item [iii)] Let $f\in L^{p}_{A}(\mathbb{R}_{+})$,  $p \in [1,+\infty]$. For all  $x\in\mathbb{R}_{+}$, the function $\tau_{x}f$
belongs to $ L^{p}_{A}(\mathbb{R}_{+})$,  $p \in [1,+\infty]$, and
we have
\begin{eqnarray*}
\|\tau_{x}f\|_{A, p}\leq\|f\|_{A,p}.
\end{eqnarray*}
\end{itemize}
\section{Generalized Fourier transform}
\label{sec:1} Throughout this section, $k$ refers to the constant
obtained in (3) and (4) from the estimates of $|c(\lambda)|^{-2}$.
\\

In the following lemma, we prove the Hardy-Littlewood inequality for
the Fourier transform.
\begin{lemma} For $f\in L^{p}_{A}(\mathbb{R}_{+}), \, 1<p\leq2$, one has
\begin{eqnarray}
\int_{\mathbb{R}_{+}}(g(x))^{p-2}|\mathcal{F}(f)(x)|^{p}\frac{dx}{|c(x)|^{2}}
\leq c\,\|f\|_{A,p}^{p}
\end{eqnarray}
where
\begin{itemize}
\item [i)]
 $g(x)=x^{2(\alpha+1)}$ if $\rho=0$ and $\alpha>0$.
\item [ii)]
$g(x)= \left\{
\begin{array}{ll}
x^{2(\alpha+1)}& \mbox{for } x>k \\
x^{3}& \mbox{for } x\leq k.
\end{array}
\right.\qquad \textrm{if}\quad \rho>0 \quad and\quad \alpha
>-\frac{1}{2}$ \\where $k$  refers to the constant
obtained from the estimates of $\,|c(x)|^{-2}.$
\end{itemize}
\end{lemma}
\begin{proof}
For $f\in L^{p}_{A}(\mathbb{R}_{+})$, $1\leq p\leq2$, we consider
the operator
$$L(f)(x)=g(x)\mathcal{F}(f)(x), \; x\in\mathbb{R}_{+}.$$
For every $f\in L^{2}_{A}(\mathbb{R}_{+})$, we have from (6)
\begin{eqnarray*}
(\int_{\mathbb{R}_{+}}|L(f)(x)|^{2}\frac{dx}{(g(x))^{2}|c(x)|^{2}})^{\frac{1}{2}}=\|\mathcal{F}(f)\|_{c,2}=\|f\|_{A,2},
\end{eqnarray*}
 hence $L$ is an operator  of strong-type $(2,2)$ between the spaces $(\mathbb{R}_{+}, A(x)dx)$
 and\\
 $(\mathbb{R}_{+}, \frac{dx}{(g(x))^{2}|c(x)|^{2}})$.
\\i) Assume $\rho=0$, $\alpha>0$ and $g(x)=x^{2(\alpha+1)}$.
For $\lambda\in]0, +\infty[$,   $f\in L^{1}_{A}(\mathbb{R}_{+})$ and
using (2) and (5), we can write
\begin{eqnarray*}
\int_{\{x\in\mathbb{R}_{+}:\,
|L(f)(x)|>\lambda\}}\frac{dx}{(g(x))^{2}|c(x)|^{2}}&=&
\int_{\{x\in\mathbb{R}_{+}:\,
|L(f)(x)|>\lambda\}}\frac{dx}{x^{4(\alpha+1)}|c(x)|^{2}}\nonumber\\&\leq&
c\int^{+\infty}_{(\frac{\lambda}{\|f\|_{A,1}})^\frac{1}{2(\alpha+1)}}\frac{x^{2\alpha+1}}{x^{4(\alpha+1)}}dx\nonumber\\&\leq&
c\,\frac{\|f\|_{A,1}}{\lambda}.
\end{eqnarray*}
It yields that  $L$ is of weak-type $(1,1)$ between the spaces under
consideration.\\ By the Marcinkiewicz interpolation theorem (see
[10]), we can assert that $L$ is an operator of strong-type $(p,p)$,
$1<p\leq2$ between the spaces $(\mathbb{R}_{+}, A(x)dx)$ and
$(\mathbb{R}_{+}, \frac{dx}{(g(x))^{2}|c(x)|^{2}})$.
\\We conclude that,
\begin{eqnarray*}
\int_{\mathbb{R}_{+}}
|L(f)(x)|^{p}\frac{dx}{(g(x))^{2}|c(x)|^{2}}&=&
\int_{\mathbb{R}_{+}}|g(x)|^{p-2}|\mathcal{F}(f)(x)|^{p}\frac{dx}{|c(x)|^{2}}\\&\leq&
c\,\|f\|_{A,p}^{p}\;\;,
\end{eqnarray*}
witch prove the result.
\\ii) Suppose now  $\rho>0$,   $\alpha >-\frac{1}{2}$ and $g(x)=
\left\{
\begin{array}{ll}
x^{2(\alpha+1)}& \mbox{for } x>k \\
x^{3}& \mbox{for } x\leq k,
\end{array}
\right.$
\\where $k$ is the constant obtained in (3) and (4) from the estimates of $|c(\lambda)|^{-2}$.
Let $\lambda\in]0, +\infty[$ and  $f\in L^{1}_{A}(\mathbb{R}_{+})$,
by (3), (4) and (5), we have
\begin{alignat*}{2}
&\int_{\{x\in\mathbb{R}_{+}:\,
|L(f)(x)|>\lambda\}}\frac{dx}{(g(x))^{2}|c(x)|^{2}}\leq
\int_{\{x\in\mathbb{R}_{+}:\,
g(x)>\frac{\lambda}{\|f\|_{A,1}}\}}\frac{dx}{(g(x))^{2}|c(x)|^{2}}\\&\leq
\int_{\{x\in\mathbb{R}_{+}:\,
g(x)>\frac{\lambda}{\|f\|_{A,1}}\}}\chi_{[0,
k]}(x)\frac{dx}{(g(x))^{2}|c(x)|^{2}}\\&+
\int_{\{x\in\mathbb{R}_{+}:\,
g(x)>\frac{\lambda}{\|f\|_{A,1}}\}}\chi_{]k,
+\infty[}(x)\frac{dx}{(g(x))^{2}|c(x)|^{2}}\\&\leq
c\int^{+\infty}_{(\frac{\lambda}{\|f\|_{A,1}})^{\frac{1}{3}}}\chi_{[0,
k]}(x)\frac{x^{2}}{x^{6}}dx+
c\int^{+\infty}_{(\frac{\lambda}{\|f\|_{A,1}})^{\frac{1}{2(\alpha+1)}}}\chi_{]k,
+\infty[}(x)\frac{x^{2\alpha+1}}{x^{4(\alpha+1)}}dx\\&\leq
c\int^{+\infty}_{(\frac{\lambda}{\|f\|_{A,1}})^{\frac{1}{3}}}x^{-4}dx+
c\int^{+\infty}_{(\frac{\lambda}{\|f\|_{A,1}})^{\frac{1}{2(\alpha+1)}}}x^{-2\alpha-3}dx\leq
c\,\frac{\|f\|_{A,1}}{\lambda}.
\end{alignat*}
Hence $L$ is of weak-type $(1,1)$ between the spaces
$(\mathbb{R}_{+}, A(x)dx)$ and $(\mathbb{R}_{+},
\frac{dx}{(g(x))^{2}|c(x)|^{2}})$.
\\We conclude by the Marcinkiewicz interpolation theorem that $L$ is of strong-type $(p,p)$, between the spaces under consideration.
\\It yields, that
\begin{eqnarray*}
\int_{\mathbb{R}_{+}}
|L(f)(x)|^{p}\frac{dx}{(g(x))^{2}|c(x)|^{2}}&=&
\int_{\mathbb{R}_{+}}|g(x)|^{p-2}|\mathcal{F}(f)(x)|^{p}\frac{dx}{|c(x)|^{2}}\\&\leq&
c\,\|f\|_{A,p}^{p}\;\;,
\end{eqnarray*}
thus we obtain the result.
\end{proof}

In the following, we study the integrability of the generalized
Fourier transform in the Jacobi hypergroup case (see Remarks 1 and
2). For $1\leq p\leq2$, we denote by $p'$ the conjugate of $p$.
\begin{lemma}
Let $1\leq p\leq2$ and $f\in L^{p}_{A}(\mathbb{R}_{+})$. Then there
exists a positive constant $c$ such that for $\delta>0$, one has
$$\Big(\int_{0}^{+\infty}\min\{1,(\delta x)^{2p'}\}\;|\mathcal{F}(f)(x)|^{p'}\frac{dx}{|c(x)|^{2}}\Big)^{\frac{1}{p'}}
\leq c\,\omega_{A,p}(f)(\delta),\; \mbox{if}\;\; 1< p\leq2$$ and
$$ess\sup_{x>0}\Big(\min\{1,(\delta x)^{2}\}\;|\mathcal{F}(f)(x)|\Big)\leq c\,\omega_{A,1}(f)(\delta),\; \mbox{if}\;\; p=1.$$
\end{lemma}
\begin{proof}
For $f\in L^{p}_{A}(\mathbb{R}_{+})$, $1\leq p\leq2$, we have by (8)
\begin{eqnarray*}
\mathcal{F}(\tau_{\delta}(f)-f)(x)=(\varphi_{x}(\delta)-1)\mathcal{F}(f)(x),
\end{eqnarray*}
for $\delta>0$ and a.e $x\in\mathbb{R}_{+}$. Applying (7), we get
\begin{eqnarray*}
\|\mathcal{F}(\tau_{\delta}(f)-f)\|_{c,p'}&=&\Big(\int_{0}^{+\infty}|1-\varphi_{x}(\delta)|^{p'}
|\mathcal{F}(f)(x)|^{p'}\frac{dx}{|c(x)|^{2}}\Big)^{\frac{1}{p'}}\\&\leq&
c\,\omega_{A,p}(f)(\delta).
\end{eqnarray*}
From (1), we obtain our results. Here  when $p=1$,  we make the
usual modification.
\end{proof}
\begin{remark}
\noindent\begin{itemize}
\item [i)] In the lemma 2, the gauge on the size of the generalized transform in terms of an integral modulus of continuity of $f$
gives a quantitative form of the Riemann-Lebesgue lemma:
$$\Big(\int_{\frac{1}{\delta}}^{+\infty}|\mathcal{F}(f)(x)|^{p'}\frac{dx}{|c(x)|^{2}}\Big)^{\frac{1}{p'}}
\leq c\,\omega_{A,p}(f)(\delta),\; \mbox{if}\;\; 1< p\leq2$$ and
$$ess\sup_{x>\frac{1}{\delta}}|\mathcal{F}(f)(x)|\leq c\,\omega_{A,1}(f)(\delta),\; \mbox{if}\;\; p=1.$$
\item [ii)] We will use the following estimates deduced from lemma 2 to establish the integrability of $\mathcal{F}(f)$ when $f$ belongs in
$\mathcal{B}^{p,\infty}_{\gamma,\alpha}$ for $1\leq p\leq2$:
\begin{eqnarray}
\delta^{2}\Big(\int_{0}^{\frac{1}{\delta}}x^{2p'}|\mathcal{F}(f)(x)|^{p'}\frac{dx}{|c(x)|^{2}}\Big)^{\frac{1}{p'}}
\leq c\,\omega_{A,p}(f)(\delta),\; \mbox{if}\; \;1< p\leq2
\end{eqnarray} and
\begin{eqnarray}ess\sup_{0<x<\frac{1}{\delta}}\Big((\delta x)^{2}|\mathcal{F}(f)(x)|\Big)\leq c\,\omega_{A,1}(f)(\delta),\; \mbox{if}\;\; p=1.\end{eqnarray}
\end{itemize}
\end{remark}
\begin{theorem}
If $f\in
\mathcal{B}^{p,1}_{\frac{2(\alpha+1)}{p},\alpha}\cap\mathcal{B}^{p,1}_{\frac{3}{p},\alpha}$
for $1<p\leq2$, then
$$\mathcal{F}(f)\in L^{1}_{c}(\mathbb{R}_{+}).$$
\end{theorem}
\begin{proof}
For $f\in L^{p}_{A}(\mathbb{R}_{+})$, $1<p\leq2$ and $\delta>0$, we
can write from (8) and (9)
\begin{eqnarray*}
\int_{\mathbb{R}_{+}}|1-
\varphi_{t}(\delta)|^{p}|\mathcal{F}(\tau_{\delta}(f)(t)|^{p}(g(t))^{p-2}\frac{dt}{|c(t)|^{2}}
&\leq& c\,(\omega_{A,p}(f)(\delta))^{p},
\end{eqnarray*}
then by (1), we obtain
\begin{eqnarray}
\delta^{2p}\int_{0}^{\frac{1}{\delta}}t^{2p}|\mathcal{F}(f)(t)|^{p}(g(t))^{p-2}\frac{dt}{|c(t)|^{2}}
\leq c\,(\omega_{A,p}(f)(\delta))^{p}.
\end{eqnarray}
From (3) and (4), we have
\\$\displaystyle\int_{0}^{\frac{1}{\delta}}t|\mathcal{F}(f)(t)|\frac{dt}{|c(t)|^{2}}$
\begin{eqnarray*}
&=&\int_{0}^{\frac{1}{\delta}}t|\mathcal{F}(f)(t)|\chi_{[0,
k]}(t)\frac{dt}{|c(t)|^{2}}+
\int_{0}^{\frac{1}{\delta}}t|\mathcal{F}(f)(t)|\chi_{]k,
+\infty[}(t)\frac{dt}{|c(t)|^{2}}\\&\leq&
c\,\int_{0}^{\frac{1}{\delta}}t|\mathcal{F}(f)(t)|\chi_{[0,
k]}(t)[t^{2}dt]+
c\,\int_{0}^{\frac{1}{\delta}}t|\mathcal{F}(f)(t)|\chi_{]k,
+\infty[}(t)[t^{2\alpha+1}dt],
\end{eqnarray*}
by H\"{o}lder's inequality and (12), we have
\\$\displaystyle\int_{0}^{\frac{1}{\delta}}t|\mathcal{F}(f)(t)|\frac{dt}{|c(t)|^{2}}$
\begin{eqnarray*}
&\leq&c\,\Big(\int_{0}^{\frac{1}{\delta}}t^{3(p-2)+2p}|\mathcal{F}(f)(t)|^{p}\chi_{[0,
k]}(t)[t^{2}dt]\Big)^{\frac{1}{p}}
\Big(\int_{0}^{\frac{1}{\delta}}t^{2(p'-2)}\chi_{[0,
k]}(t)dt\Big)^{\frac{1}{p'}}\\&+&
c\,\Big(\int_{0}^{\frac{1}{\delta}}t^{2(\alpha+1)(p-2)+2p}|
\mathcal{F}(f)(t)|^{p}\chi_{]k,
+\infty[}(t)[t^{2\alpha+1}dt]\Big)^{\frac{1}{p}}
\\&\times&\Big(\int_{0}^{\frac{1}{\delta}}t^{(2\alpha+1)(p'-2)+2\alpha-1}\chi_{]k, +\infty[}(t)dt\Big)^{\frac{1}{p'}}\\&\leq&
c\,\Big(\int_{0}^{\frac{1}{\delta}}t^{2p}|\mathcal{F}(f)(t)|^{p}(g(t))^{p-2}\frac{dt}{|c(t)|^{2}}\Big)^{\frac{1}{p}}
\\&\times&\Big\{\Big(\int_{0}^{\frac{1}{\delta}}t^{2(p'-2)}dt\Big)^{\frac{1}{p'}}+
\Big(\int_{0}^{\frac{1}{\delta}}t^{(2\alpha+1)(p'-2)+2\alpha-1}dt\Big)^{\frac{1}{p'}}\Big\}\\&\leq&
c\,\delta^{-2}\omega_{A,p}(f)(\delta)\Big(\frac{1}{\delta^{\frac{3}{p}-1}}+
\frac{1}{\delta^{\frac{2(\alpha+1)}{p}-1}}\Big)\leq
c\,\Big(\frac{\omega_{A,p}(f)(\delta)}{\delta^{\frac{3}{p}}}\frac{1}{\delta}+
\frac{\omega_{A,p}(f)(\delta)}{\delta^{\frac{2(\alpha+1)}{p}}}\frac{1}{\delta}\Big).
\end{eqnarray*}
Integrating with respect to $\delta$ over $\mathbb{R}_{+}$ for $f\in
\mathcal{B}^{p,1}_{\frac{2(\alpha+1)}{p},\alpha}\cap
\mathcal{B}^{p,1}_{\frac{3}{p},\alpha}$, the double integral is
evaluated by interchanging the orders of integration, it yields
$$\int_{0}^{+\infty}|\mathcal{F}(f)(t)|\frac{dt}{|c(t)|^{2}}<+\infty.$$
This complete the proof.
\end{proof}
\begin{theorem}
Let $\gamma>0$, $1\leq p \leq 2$ and  $f\in
\mathcal{B}^{p,\infty}_{\gamma,\alpha}$, then
\begin{itemize}
\item [i)] For $p\neq1$ and $0<  \gamma \leq \frac{2(\alpha+1)}{p}$,
one has
\\$\mathcal{F}(f)\in  L^{s}_{c}(\mathbb{R}_{+})$ provided that
$\frac{2(\alpha+1)p}{\gamma p+2(\alpha+1)(p-1)}< s \leq p'.$
\item [ii)] For $p\neq1$ and $\gamma >\frac{2(\alpha+1)}{p}$, one
has
\begin{eqnarray*}
\mathcal{F}(f)\in  L^{1}_{c}(\mathbb{R}_{+}).
\end{eqnarray*}
\item [iii)] For $p=1$ and  $\gamma>\sup(3,2(\alpha+1))$, one has
\begin{eqnarray*}
\mathcal{F}(f)\in  L^{1}_{c}(\mathbb{R}_{+}).
\end{eqnarray*}
\end{itemize}
\end{theorem}
\begin{proof} Let $f\in \mathcal{B}^{p,\infty}_{\gamma,\alpha}$, $1\leq p \leq 2$.
\\i) Suppose that  $p\neq1$ and $0<  \gamma \leq \frac{2(\alpha+1)}{p}$.
Let $\frac{2(\alpha+1)p}{\gamma p+2(\alpha+1)(p-1)}< s \leq p'$, we
define the function
$$g(t)=\int_{k}^{t}|\mathcal{F}(f)(x)|^{s}x^{s}\frac{dx}{|c(x)|^{2}},\quad t>k.$$
By H\"{o}lder's inequality, (4) and (10) we have
\begin{eqnarray*}
g(t)&\leq&
\Big(\int_{k}^{t}|\mathcal{F}(f)(x)|^{p'}x^{2p'}\frac{dx}{|c(x)|^{2}}\Big)^{\frac{s}{p'}}
\Big(\int_{k}^{t}\frac{dx}{|c(x)|^{2}}\Big)^{1-\frac{s}{p'}}\\&\leq&
c\,t^{2s}(\omega_{A,p}(f)(\frac{1}{t}))^{s}\Big(\int_{k}^{t}\frac{dx}{|c(x)|^{2}}\Big)^{1-\frac{s}{p'}}\\&\leq&
c\,t^{(2-\gamma)s}\Big(\int_{k}^{t}x^{2\alpha+1}dx\Big)^{1-\frac{s}{p'}}\leq
c\,t^{(2-\gamma)s+2(\alpha+1)(1-\frac{s}{p'})}.
\end{eqnarray*}
Then we get
\begin{eqnarray*}
\int_{k}^{t}|\mathcal{F}(f)(x)|^{s}\frac{dx}{|c(x)|^{2}}&=&\int_{k}^{t}x^{-2s}g'(x)dx\\&=&
t^{-2s}g(t)+2s\int_{k}^{t}x^{-2s-1}g(x)dx\\&\leq& c\,\Big(t^{-\gamma
s+2(\alpha+1)(1-\frac{s}{p'})}+\int_{k}^{t}x^{-\gamma
s+2(\alpha+1)(1-\frac{s}{p'})-1}dx\Big)\\&\leq& c\,(t^{-\gamma
s+2(\alpha+1)(1-\frac{s}{p'})}+1),
\end{eqnarray*}
it yields that $\mathcal{F}(f)\in  L^{s}_{c}(]k,+\infty[,
\frac{dx}{|c(x)|^{2}})$. Since $\mathcal{F}(f)\in L^{p'}([0,k],
\frac{dx}{|c(x)|^{2}})\subset L^{s}([0,k], \frac{dx}{|c(x)|^{2}})$,
we deduce that $\mathcal{F}(f)$ is in $L^{s}_{c}(\mathbb{R}_{+})$.
\\ii) Assume now $\gamma >\frac{2(\alpha+1)}{p}$. For $p\neq1$, by proceeding in the same manner as the proof of i)
with $s=1$, we obtain the desired result.
\\iii) For $p=1$ and $\gamma>\sup(3,2(\alpha+1))$. By H\"{o}lder's inequality, (3), (4) and (11), we
have for $t>0$
\begin{eqnarray*}
\int_{0}^{\frac{1}{t}}|\mathcal{F}(f)(x)|x\frac{dx}{|c(x)|^{2}}&\leq&
ess\sup_{0<x\leq\frac{1}{t}} x^{2}|\mathcal{F}_k(f)(x)|
\int_{0}^{\frac{1}{t}}\frac{1}{x}\frac{dx}{|c(x)|^{2}}\\&\leq&
c\,t^{\gamma-2}\Big(\int_{0}^{\frac{1}{t}}\frac{1}{x}\chi_{\{0\leq
x\leq k \}}\frac{dx}{|c(x)|^{2}}+
\int_{0}^{\frac{1}{t}}\frac{1}{x}\chi_{\{ x> k
\}}\frac{dx}{|c(x)|^{2}}\Big)\\&\leq&
c\,t^{\gamma-2}[t^{-2}+t^{-(2\alpha+1)}]\leq
c\,[t^{(\gamma-3)-1}+t^{\gamma-2(\alpha+1)-1}].
\end{eqnarray*}
Integration with respect to  $t$ over $ (0,1)$  and applying
Fubini's theorem we obtain
$$\int_{1}^{+\infty}|\mathcal{F}(f)(x)|\frac{dx}{|c(x)|^{2}}\leq
c\Big(\int_{0}^{1}t^{(\gamma-3)-1}dt+\int_{0}^{1}t^{\gamma-2(\alpha+1)-1}dt\Big)<\infty.$$
Since $L^{\infty}([0,1], \frac{dx}{|c(x)|^{2}})\subset
L^{1}([0,1],\frac{dx}{|c(x)|^{2}})$, then $\mathcal{F}(f)\in
L^{1}_{c}(\mathbb{R}_{+})$.
\end{proof}
\begin{remark} For $\gamma>\sup(3,2(\alpha+1))$, we can assert from the theorem 2, iii)
that $\mathcal{B}^{1,\infty}_{\gamma,\alpha}$ is an example of space
where we can apply the inversion formula.
\end{remark}
\begin{acknowledgements}
The authors are supported by the DGRST research project 04/UR/15-02
and the program CMCU 10G / 1503.
\end{acknowledgements}

\bibliographystyle{spmpsci}      


\end{document}